\documentclass[11pt,twoside]{amsart}
\usepackage{amsmath}
\usepackage{amssymb}
\usepackage{latexsym}
\usepackage{graphics}
\usepackage[all]{xy}
\usepackage{url} 
\usepackage{graphicx}
\usepackage{float}

\newtheorem{thm}{Theorem}

\newtheorem{df}{Definition}

\begin{document}

\leftline{ \scriptsize}

\vspace{1.3 cm}
\title[Inclusions of  Waterman-Shiba  spaces into  generalized Wiener classes]
{Inclusions of  Waterman-Shiba  spaces\\
 into  generalized Wiener classes}
\author{Mahdi Hormozi}
\address{Department of Mathematical Sciences, Division of Mathematics\\ University of Gothenburg\\ Gothenburg 41296\\ Sweden}
\email{hormozi@chalmers.se}
\author{Franciszek Prus-Wi\'{s}niowski}
\address{Institute of Mathematics\\ Szczecin University\\ ul. Wielkopolska 15\\ PL-70-451 Szczecin\\ Poland}
\email{wisniows@univ.szczecin.pl}
\author{Hjalmar Rosengren}
\address{Department of Mathematical Sciences, Division of Mathematics\\  Chalmers University of Technology\\ Gothenburg 41296\\ Sweden}
\email{hjalmar@chalmers.se}
\thanks{{\scriptsize
\hskip -0.4 true cm MSC(2000):Primary  26A15; Secondary 26A45
}}
\hskip -0.4 true cm

\maketitle



\begin{abstract} The characterization of the inclusion of  Waterman-Shiba spaces $\:\Lambda BV^{(p)}\:$ into generalized Wiener
classes of functions  $BV(q;\,\delta)$ is given. It  uses a new  and shorter proof and extends an earlier result of U. Goginava.

\end{abstract}

\vskip 0.2 true cm


%


\vskip 0.4 true cm
Let $\:\Lambda=(\lambda_i)\:$ be a $\Lambda$-sequence, that is, a nondecreasing sequence of positive numbers such that
$\:\sum\frac1{\lambda_i}=+\infty\:$ and let $p$ be a number greater than or equal to 1. A function $f:[0,\,1]\to\mathbb R$ is said to be of bounded $p$-$\Lambda$-variation
 if
$$
 V(f)\ :=\ \sup\left(\sum_{i=1}^n\frac{|f(I_i)|^p}{\lambda_i}
\right)^\frac1p\ <\ +\infty,
$$
where the supremum is taken over all finite families $\:\{I_i\}_{i=1}^n\:$ of
nonoverlapping subintervals of  $[0,1]$ and where $f(I_i):=f(\sup I_i)-f(\inf I_i)$
is the change of the function $f$ over the interval $I_i$. The symbol $\Lambda BV^{(p)}$ denotes
the linear space of all functions of bounded $p$-$\Lambda$-variation with domain $[0,\,1]$. The Waterman-Shiba space
$\:\Lambda BV^{(p)}\:$ was introduced in 1980 by M. Shiba  \cite{S}. When $p=1$, $\:\Lambda BV^{(p)}\:$ is the well-known Waterman space $\:\Lambda BV$ (see e.\ g.\ \cite{Wt1} and \cite{Wt2}). Some of the properties and applications of functions of class $\:\Lambda BV^{(p)}\:$ were discussed in \cite{BT}, \cite{BTV},  \cite{Hor},  \cite{HLP}, \cite{L}, \cite{Lind},  \cite{PW}, \cite{SW1}, \cite{SW}, \cite{V1}, \cite{V2} and \cite{V3}. $\Lambda BV^{(p)}\:$ equipped with the norm $\:\|f\|_{\Lambda,\,p}:=
|f(0)|+V(f)\:$ is a Banach space.\\

H. Kita and K. Yoneda introduced a new function space which is a generalization of Wiener classes \cite{KY} ( see also \cite{Ki} and \cite{Akh}). The concept was further extended by T. Akhobadze in \cite{Akh2} who studied many properties of the generalized Wiener classes $\:BV(q,\,\delta)$ thoroughly (see \cite{Akh3}, \cite{Akh4}, \cite{Akh5}, \cite{Akh6}, \cite{Akh7}).

\begin{df}\label{def1}
 Let $q=(q(n))_{n=1}^\infty$ be an increasing positive sequence and let $\delta=(\delta(n))_{n=1}^\infty$ be an increasing and unbounded positive sequence. We say that a function $f:[0,\,1]\to\mathbb R$  belongs to the class
$BV(q;\delta)$ if
$$
V(f,  q;\, \delta):=\sup_{n\geq 1}\sup_{\{I_k\}} \left\{ \left(\sum_{k=1}^{s} |f(I_k)|^{q(n)}\right)^{\tfrac{1}{q(n)}} :\inf_{k} |I_k|\geq \frac{1}{\delta(n)} \right\} < \infty,
$$
where $\{I_k\}_{k=1}^s$ are non-overlapping subintervals
of $[0,1]$.
\end{df}

If $\delta(n)^{1/q(n)}$ is a bounded sequence,
then $BV(q;\delta)$ is simply the space of all bounded functions. This
follows from the estimate
$$\left(\sum_{k=1}^{s} |f(I_k)|^{q(n)}\right)^{\tfrac{1}{q(n)}}
\leq 2Cs^{\frac 1{q(n)}}\leq 2C\delta(n)^{^\frac 1{q(n)}},\qquad |f(x)|\leq C.$$

 The following statement regarding inclusions of Waterman spaces into  generalized Wiener classes has been presented in \cite {G3}: if $\lim_{n\rightarrow\infty}q(n)=\infty$ and $\delta(n)=2^n$ the inclusion $\:\Lambda BV\subset BV(q,\delta)\:$ holds if and only if
   \begin{equation}
    \label{e1}
    \limsup_{n\rightarrow\infty}  \left \{\ \max_{1\leq k\leq 2^n} \ \frac{k^{\frac{1}{q(n)}}}{(\sum_{i=1}^{k} \frac{1}{\lambda_i})} \right\}< +\infty.
   \end{equation}
Our result formulated below extends the above theorem of Goginava essentially and furnishes a new and much shorter proof.

    \begin{thm}
   \label{t1}
   For $p \in [1 , \infty)$ and $q$ and $\delta$  sequences satisfying
the conditions in Definition \ref{def1},
 the inclusion $\Lambda BV^{(p)}\subset BV(q;\,\delta)$ holds if and only if
   \begin{equation}
    \label{e2}
   \limsup_{n\rightarrow\infty}  \left \{\ \max_{1\leq k\leq \delta (n)} \ \frac{k^{\frac{1}{q(n)}}}{(\sum_{i=1}^{k} \frac{1}{\lambda_i})^{\frac{1}{p}}} \right\}< +\infty.
   \end{equation}
   \end{thm}

 Before we present a relatively short proof of Thm.\ \ref{t1}, we
give an example showing that it provides a non-trivial extension
of \eqref{e1} even for $p=1$. The  example \--- provided by the referee kindly instead of our more complicated one \--- is obtained by taking $\:\lambda_n=n$, $\:q(n)=\sqrt{n}\;$ and $\: \delta(n)=2^{\sqrt n}$. With those choices, it follows immediately that the Goginava indicator \eqref{e1} is infinite while \eqref{e2} holds.

\begin{proof}[Proof of Thm. \ref{t1}]
To show that  (\ref{e2}) is a \textbf{sufficiency} condition for the inclusion $\Lambda BV^{(p)}\subset BV(q;\,\delta)$, we will
prove the inequality
\begin{equation}\label{vvin} V(f, q,\delta) \leq V_{\Lambda BV^p}(f) \sup_{n}\left \{\max_{1\leq k\leq \delta(n)}\ \frac{k}{(\sum_{i=1}^{k} 1/\lambda_{i})^{\frac{q(n)}{p}}} \right\} ^{\frac{1}{q(n)}}.\end{equation}
This is a consequence of the numerical inequality
\begin{equation}\label{homin}\sum_{j=1}^sx_j^q\leq\left(\sum_{j=1}^s{x_jy_j}\right)^q
\max_{1\leq k\leq s}\frac{k}{\left(\sum_{j=1}^k{y_j}\right)^q},
 \end{equation}
which is valid for $q\geq 0$ and
\begin{align}\label{xd}x_1&\geq x_2\geq \dots\geq x_s\geq 0,\\
\notag y_1&\geq y_2\geq \dots\geq y_s\geq 0.
\end{align}
 In the cases $q\geq 1$ and $0\leq q<1$, \eqref{homin} is a
reformulation of
\cite[Lemma]{kup} and  \cite[Lemma 2.5]{Hor}, respectively.

To prove \eqref{vvin}, consider a
non-overlapping family  $(I_k)_{k=1}^s$ with
$\inf|I_k|\geq 1/\delta(n)$. In particular, $s\leq \delta(n)$.
Apply  \eqref{homin} with $q$ replaced by $q(n)/p$,
 $x_j=|f(I_j)|^p$ and $y_j=1/\lambda_j$, where we
 reorder the intervals so that \eqref{xd} holds.
This gives
\begin{align*}\left(\sum_{j=1}^s|f(I_j)|^{q(n)}\right)^{\frac 1{q(n)}}
&\leq \left(\sum_{j=1}^s\frac{|f(I_j)|^p}{\lambda_j}\right)^{\frac 1p}\max_{1\leq k\leq s}\frac{k^{\frac 1{q(n)}}}{(\sum_{j=1}^k\frac 1{\lambda_j})^{\frac 1p}}\\
&\leq V_{\Lambda BV^p}(f)\max_{1\leq k\leq \delta(n)}\frac{k^{\frac 1{q(n)}}}{(\sum_{j=1}^k\frac 1{\lambda_j})^{\frac 1p}}.\end{align*}
Taking the supremum over $n$ yields
 \eqref{vvin}.

  \textbf{ Necessity.} Suppose (\ref{e2}) doesn't hold. Then there is an increasing sequence $\:(n_k)\:$ of positive integers such that for all indices $k$
 \begin{equation}
\label{e6}
\delta(n_k) \geq 2^{k+2}
\end{equation}
and
 \begin{equation}
\label{e7}
\max_{1\le n\le \delta(n_k)}\,\frac{n^\frac1{q(n_k)}}{\left(\sum_{i=1}^n\frac1{\lambda_i}\right)^\frac1p}\ >\ 2^{2k+\frac{k+1}{q(1)}}.
\end{equation}
Let $\:(m_k)\:$ be a sequence of positive integers such that
\begin{equation}
\label{e8}
1\ \le\ m_k\ \le\ \delta(n_k),
\end{equation}
and
\begin{equation}
\label{e9}
\max_{1\le n\le \delta(n_k)}\,\frac{n^\frac1{q(n_k)}}{\left(\sum_{i=1}^n\frac1{\lambda_i}\right)^\frac1p}\ \ =\ \ \frac{m_k^\frac1{q(n_k)}}{\left(\sum_{i=1}^{m_k}\frac1{\lambda_i}\right)^\frac1p}.
\end{equation}
 Denote
$$
 \Phi_{k} :=  \frac{1}{\sum_{i=1}^{m_k} 1/\lambda_{i}}.
$$
Consider
$$ g_k(y):=\begin{cases}
 2^{-k}\Phi_{k}^{1/p} ~~~~~, ~~~~y\in [\tfrac1{2^k}+\frac{2j-2}{\delta(n_k)},\tfrac1{2^k}+\frac{2j-1}{\delta(n_k)}) ;~~~~~  1\leq j\leq N_k, \\[.1in]
0~~~~~\qquad\qquad \textmd{otherwise},\\
\end{cases}
$$
where
$$N_k= \min \{m_k, s_k\},\qquad
s_k\ =\ \max \left\{j\in \mathbb{N}:\ 2j \leq \frac{\delta(n_k)}{2^k}+1\,\right\}.
$$
Applying the fact that $2(s_k+1)\geq \frac{\delta(n_k)}{2^k}+1$ and (\ref{e6}), we have
\begin{equation}
\label{e10}
 \frac{2s_k-1}{\delta(n_k)} \geq 2^{-k-1}.
\end{equation}

The functions $g_k$ have disjoint support and the series $\:\sum_{k=1}^{\infty}g_k(x)$ converges uniformly to a function $g$. Thus,
\begin{align*}
\|g\|_{\Lambda,p} & \le\ \sum_{k=1}^{\infty}\|g_k\|_{\Lambda,p}\
= \sum_{k=1}^{\infty}  \left(\sum_{j=1}^{2N_k} \frac{(2^{-k}\Phi^{1/p}_{k})^p}{\lambda_j}\right)^{1/p}\leq \sum_{k=1}^{\infty} 2^{-k} \left(2 \sum_{j=1}^{N_k} \frac{\Phi_{k}}{\lambda_j}\right)^{1/p} \\[.1in]
&
\leq\ \sum_{k=1}^{\infty} 2^{-k} \left(2 \sum_{j=1}^{m_k} \frac{\Phi_{k}}{\lambda_j}\right)^{1/p} \
= \ \sum_{k=1}^{\infty} 2^{-k} \left(2\  \frac{ \sum_{j=1}^{m_k} \frac{1}{\lambda_j}}{ \sum_{j=1}^{m_k} \frac{1}{\lambda_j}}\right)^{1/p}
  \ <\ +\infty,
\end{align*}
that is,  $g \in \Lambda BV^{(p)}$.

If $N_k=m_k$, then $\:2N_k-1\ge\ m_k$, and if $\:N_k=s_k$, then
$$
2N_k-1\ \overset{\text{\eqref{e10}}}{\ge}\ \ \frac{\delta(n_k)}{2^{k+1}}\ \ \overset{\text{\eqref{e8}}}{\ge} \ \frac{m_k}{2^{k+1}}.
$$
Hence
\begin{equation}
\label{e11}
2N_k-1\ \ge\ \frac{m_k}{2^{k+1}} \hspace{.3in}\text{for all $k$.}
\end{equation}

 Now, given any positive integer $k$, all intervals  $[\tfrac1{2^k}+\frac{j-1}{\delta(n_k)}, \tfrac1{2^k}+\frac{j}{\delta(n_k)}]$ for $j=1,..,2N_k-1$, have  length  $\frac{1}{\delta(n_k)}\:$ and thus
\begin{align*}
V(g, q,\delta)  &\ge \  \left( \sum_{j=1}^{2N_k-1} \bigl|g\left(\tfrac1{2^k}+\tfrac{j-1}{\delta(n_k)}\right)-g\left(\tfrac1{2^k}+\tfrac{j}{\delta(n_k)}\right)\bigr|^{q(n_k)} \right)^{\frac1{q(n_k)}}\\[.1in] &=\ \left( (2N_k-1)(2^{-k}\Phi_k^\frac1p)^{q(n_k)}\right)^{\frac1{q(n_k)}}\ =\ \frac1{2^k}\left(\frac{2N_k-1}{\left(\sum_{i=1}^{m_k}\frac1{\lambda_i}\right)^{\frac{q(n_k)}p}}\right)^{\frac1{q(n_k)}}\\[.1in]
&\overset{\text{\eqref{e11}}}{\ge}\ \frac1{2^k}\left(\frac1{2^{k+1}}\cdot\frac{m_k}{\left(\sum_{i=1}^{m_k}\frac1{\lambda_i}\right)^{\frac{q(n_k)}p}}\right)^{\frac1{q(n_k)}}\overset{\text{\eqref{e7}, \eqref{e9}}}{\ge}\ \frac{2^{k+(k+1)/q(1)}}{2^{(k+1)/q(n_k)}}
{\ge}\ 2^k.
\end{align*}
Since $k$ was arbitrary, $V(g, q,\delta)$ must be infinite which shows that $g\notin BV(q;\,\delta)$.
\end{proof}

\textbf{Acknowledgment.}  Hjalmar Rosengren is supported by the Swedish Science Research Council.

\end{document}